\documentclass[11pt]{article}

\usepackage[latin1]{inputenc}
\usepackage{latexsym}
\usepackage{amsfonts}
\usepackage{amsmath}
\usepackage{graphicx}

\newtheorem{theorem}{Theorem}

\newtheorem{lemma}[theorem]{Lemma}
\newtheorem{definition}[theorem]{Definition}
\numberwithin{theorem}{section}

\def\bs1{\boldsymbol\Sigma_1}

\def\et{\eta}
\def\D{\Delta}
\def \G{\Gamma }

\def\P{\mathbb{P}}
\def\F{\mathbb{F}}

\def\l{\ell}
\def\u{\operatorname{U}}

\def\Z{\mathbb{Z}}

\def\F{\mathbb{F}}
\def\C{\mathbb{C}}
\def\N{\mathbb{N}}
\def\R{\mathbb{R}}

\def\mz{\mathbf{m}_0}
\def\mu{\mathbf{m}_1}
\def\ml{\mathbf{m}_{\ell}}
\def\t{\mathbf{t}}
\def\la{\lambda}
\def\Ac{\mathcal{A}}
\def\Bc{\mathcal{B}}

\def\x{\mathbf{x}}
\def\z{\mathbf{z}}
\def\u{\mathbf{u}}
\def\v{\mathbf{v}}
\def\w{\mathbf{w}}

\def\m{\mathbf{m}}

\begin{document}

\title{On Some Explicit Semi-stable Degenerations of Toric Varieties}
\author{Marina Marchisio, Vittorio Perduca}
\date{}
\maketitle

\begin{abstract}
We study semi-stable degenerations of toric varieties determined by
certain partitions of their moment polytopes. Analyzing their defining equations we prove a property of uniqueness.

\end{abstract}

\section{Background}

\subsection{\bf Polytopes and semi-stable partitions}
 In his paper \cite{Hu}, Hu provides a toric
construction for semi-stable degenerations of toric varieties. We
study the uniqueness of this construction for a toric variety $X$ in
the particular case of a semi-stable partition of its moment polytope
in two subpolytopes. Adapting a theorem by Strumfels on toric ideals
(Lemma 4.1 in \cite{StrumL} and Section 2 in \cite{Sottile}) to
particular open polytopes, we investigate the equations of the
degeneration of $X$ as embedded variety.

\bigskip

Let   $M \simeq \Z ^n$  be a lattice and $N$ its  dual. We
consider polytopes $\D \subset M$ which describe smooth algebraic
varieties $X_\D$; $\D$ determines the normal fan $\Sigma_{X_\D}
\subset N$. Recall that convex polytopes $\D$ determine a toric
manifold $X_\D$ together with an ample line bundle $\mathcal{L}_\D$: $(X_\D,\mathcal{L}_\D)$.
If the polytope is non singular of dimension $n$, then
$\mathcal{L}_\D$ is very ample, we then have an embedding $X_\D
\hookrightarrow \P^{\l}$, for some $\l$ \cite{Oda}.

Now fix a (compact) polytope $\D$ and suppose $\D\cap M=\{\mz,\ldots,\ml\}$. Take $x_0,\ldots,x_l$ as homogeneous coordinates in $\P^\l$.
 We can define $X=X_\D$ as the closure in $\P^{\l}$ of the image of the map
\begin{eqnarray}
\varphi: (\C^*)^n & \rightarrow & \P^{\l} \label{toric} \\
          \t & \mapsto & [\t^{\mz},\ldots,\t^{\ml}], \nonumber
\end{eqnarray}
where $\t=(t_1,\ldots,t_n)\in(\C^*)^n$ and given $\u=(u_1,\ldots,u_n)\in\Z^n$ we use the notation $\t^\u=t_1^{u_1}\cdot\ldots\cdot t_n^{u_n}$. Taking homogeneous coordinates in $X_{\D}$, this map extends to a map $X_{\D}\rightarrow \P^{\l}$, which is an embedding under the assumption $X_{\D}$ smooth (see \cite{Cox}). 

\bigskip

We assume that there exists
 a suitable finite partition $\G$
of $\Delta$ in subpolytopes $\{\D_j\}_{j=1} ^k$. We will assume that
the toric varieties $X_{\D_j}$ corresponding to each $\D_j$ are also
smooth. We call an open $l$-face $\sigma$ of $\D_j$ an $l$-face of $\G$ and we declare that the $0$-faces of $\D$ are \textit{not} $0$-faces of $\G$. Following \cite{A, Hu} we ask $\G$ to be {\it semi-stable}:

\begin{definition}
$\G$ is semi-stable if for any $l$-face $\sigma$ of $\G$, if $\theta$ is a $k$-face of $\D$ such that $\sigma\subset\theta$, then there are exactly $k-l+1$  $\D_j$'s such that $\theta$ is a face of each of them.
\end{definition}

In fact:

\begin{theorem}\label{mah}\emph{\cite{A, Hu}} If $\{ \D_j \}_{j=1}^k$ is a  semi-stable partition of
$\D$, then there exists a semi-stable degeneration of $X$, $f: \tilde
X \to \C$ with central fiber $f^{-1}(0)=\cup _{j=1}^{k} X_{\D_j}$; the central
fiber is completely described by the polytope partition $\{ \D_j
\}_{j=1}^k$.
\end{theorem}

$\tilde X$ is constructed by a {\it lift} of $\D$ (see Definition (\ref{lifting})).
From Theorem 2.8 in \cite{Hu}, $\tilde X$ is unique: we study
the uniqueness of $\tilde X$ for semi-stable partitions of $\D$ in
two subpolytopes $\D_1,\D_2$, and we describe its defining
equations.
 In particular, in Section 2 of \cite{Hu}, Hu shows that
the ordering (arbitrarily fixed) $\{ \D _1, \ldots, \D_k\}$ of the
polytopes in $\G$ determines a piecewise affine function on the
partition $F: \D \to \R$, which takes rational values on the points
in the lattice $M$. $F$ can be chosen to be concave and it is called
{\it lifting function}. 

\begin{definition}
\label{lifting}
$$
\tilde \D_F= \{ ( m, \tilde m) \in M \times \Z \mbox{ such that } m \in \Delta \text{ and } \tilde m \geq F(m) \}
$$
is an {\it open lifting} (here simply {\it lift})  of $\D$ with respect
to $\G$. 
\end{definition}

There are many possible lifts of $\D$ with respect to $\G$;
if $\G$ consists of two subpolytopes, then two lifts exist. By
construction there exists a morphism $f: \tilde X_F:=X_{\tilde\D_F} \to \C$ which
realizes a semi-stable degeneration of $X$. As before we have
embeddings $X \hookrightarrow \P^{\l}$ and $\tilde X_F \hookrightarrow \P^\l \times \C$.
In particular we can define $\tilde X_F$ as the closure in $\P^{\l}\times\C$ of the image of the map:
\begin{eqnarray}
\psi_F: (\C^*)^n\times\C & \rightarrow & \P^{\l}\times\C \label{par} \\
        (\t,\la) & \mapsto & ([\la^{F(\mz)}\t^{\mz},\la^{F(\mu)}\t^{\mu},\ldots,\la^{F(\ml)}\t^{\ml}],\la). \nonumber 
\end{eqnarray} Theorem 2.8 in \cite{Hu}
claims that the image of $\psi:=\psi_F$, and hence $\tilde X_F$, is
independent of the lifting function $F$.

\bigskip

We explicitly study this statement  for semi-stable partitions of
$\D$ in two subpolytopes. If $\Gamma$ consists of two subpolytopes $\D_1,\D_2$, then we can construct two possible lifting functions $F,G$ and then $\D$ has two lifts, say $\tilde\D_F$ and $\tilde\D_G$. In particular let $y_1,\ldots,y_n$ be coordinates in $\R^n\supset\D$ and let 
$$
a_1y_1+\ldots+a_ny_n+a_{n+1}=0
$$
be an equation of the {\it cut} $\D_1\cap\D_2$ {\it in the lattice}, where we take $a_1,\ldots,a_{n+1}\in\Z$ such that for all $\m_j=(m_{1j},\ldots,m_{nj})\in \D_2\cap M$ we have 
$$
a_1m_{1j}+\ldots+a_nm_{nj}+a_{n+1}\geq 0.
$$

\bigskip

 Following the construction in \cite{Hu}, the functions $F,G$ we obtain look like:

$$
F(\m_j)=
\left\{
\begin{array}{lc}
0 & \mbox{if } \m_j\in\D_1 \\
L_F(\m_j):=a_1m_{1j}+\ldots+a_nm_{nj}+a_{n+1} & \mbox{if } \m_j\in\D_2,
\end{array}
\right. 
$$
$$
G(\m_j)=
\left\{
\begin{array}{lc}
L_G(\m_j):=-a_1m_{1j}-\ldots -a_nm_{nj}-a_{n+1} & \mbox{if } \m_j\in\D_1 \\
0 & \mbox{if } \m_j\in\D_2.
\end{array}
\right.  
$$

We prove that the two non-compact toric varieties defined by the open polytopes $\tilde\D_F$ and $\tilde\D_G$ have the same toric ideals. To do this we adapt a
Strumfels's theorem on toric ideals (Lemma 4.1 in \cite{StrumL} and
Section 2 in \cite{Sottile}) to this non-compact context.

\subsection{\bf Toric ideals}

In \cite{Sottile} Sottile describes the ideal $I$ of the compact toric variety $X$ ({\it toric ideal}) defined as the  closure of the image of a map (\ref{toric}), following Strumfels's book \cite{StrumL}.

Take $x_0,\ldots,x_l$ as homogeneous coordinates in $\P^\l$. With the notation of the previous section, suppose $\mathbf{m}_j=(m_{1j},\ldots,m_{nj})$, $j=0,\ldots,\l$ and consider the $(n+1)\times (\l+1)$ matrix
$$
\Ac^+=
\left(
\begin{array}{cccc}
1        & 1      & \ldots & 1         \\
m_{10}   & m_{11} & \ldots & m_{1\l}   \\
\vdots   & \vdots &        & \vdots    \\
m_{n0}   & m_{n1} & \ldots & m_{n\l}
\end{array}
\right).
$$

Observe that if $\u\in\Z^{\l+1}$, then we may write $\u$ uniquely as $\u=\u^+-\u^-$, where $\u^+,\u^-\in\N^{\l+1}$, but $\u^+$ and $\u^-$ have no non-zero components in common. For instance, if $\u=(1,-2,1,0)$, then $\u^+=(1,0,1,0)$ and $\u^-=(0,2,0,0)$ (Sottile's notation).

We therefore have:
\begin{theorem}\label{idc}\emph{(\cite{Sottile}, Corollary 2.3)}
$$
I= \langle \x^{\u^+}-\x^{\u^-}| \u\in\ker(\Ac^+) \mbox{ \emph{and} } \u\in\Z^{\l+1}\rangle.
$$
\end{theorem}

There are no simple formulas for a finite set of generators of a general toric ideal. An effective method for computing a finite set of equations defining $X_{\D}$ in $\P^{\l}$ is applying elimination theory algorithms to its parametrization in homogeneous coordinates. These algorithms are implemented in the well known computer algebra system Maplesoft \cite{Maple}.

\section{First examples}

To illustrate the previous section, we describe the semi-stable
degenerations of a curve and a surface determined by a
subdivision of their moment polytopes in two subpolytopes. 

\subsection{\bf The twisted cubic}

The twisted cubic $X\subset\P^3$ can be defined as $\P^1$ embedded in $\P^3$ by cubics, that is, as the toric
 curve $(X_\D,\mathcal{L}_\D)=(\P^1,\mathcal{O}(3))$, where $\D$ is the polytope below.

\begin{figure}[h]
     \begin{center}
     \scalebox{1}{\includegraphics{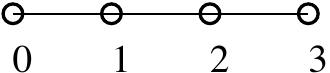}}
\caption{The moment polytope $\D$ of the twisted cubic $X\subset\P^3$.}\label{TC}
     \end{center}
     \end{figure}

Here $M=\Z$, $\D\cap M=\{\m_j=j,\,j=0,\ldots,3\}$, $X$ is the closure of the image of
\begin{eqnarray}
\varphi: \C^* & \rightarrow & \P^3 \nonumber \\
          t & \mapsto & [1,t,t^2,t^3], \nonumber
\end{eqnarray}
which extends to the embedding
$$
\begin{array}{cll}
X_{\D} & \hookrightarrow & \P^3 \nonumber \\
(v_0,v_1) & \mapsto & [v_1^3,v_0v_1^2,v_0^2v_1,v_0^3], \nonumber
\end{array}
$$
where $v_0,v_1$ are homogeneous coordinates in $X_{\D}$.

The toric ideal of $X$ is of course computed to be
$$
I=\langle x_0x_2-x_1^2, x_1x_3-x_2^2, x_0x_3-x_1x_2 \rangle.
$$

Now consider the semi-stable partition $\{\D_1,\D_2\}$ of $\D$,\
where $\D_1=[0,1]\subset\R$ and $\D_2=[1,3]\subset\R$.
This partition gives the semi-stable degeneration of $X$ to the
union of two curves $X_1\cup X_2$, where $X_1=(\P^1,\mathcal{O}(1))$
and $X_2=(\P^1,\mathcal{O}(2))$.

The two possible lifting functions are
$$
F(j)=
\left\{
\begin{array}{cc}
0 & j=0,1 \\
j-1 & j=2,3
\end{array}
\right., \,
G(j)=
\left\{
\begin{array}{cc}
1 & j=0 \\
0 & j\neq 0
\end{array}
\right..
$$

\begin{figure}[htbp]
     \begin{center}
     \scalebox{1}{\includegraphics{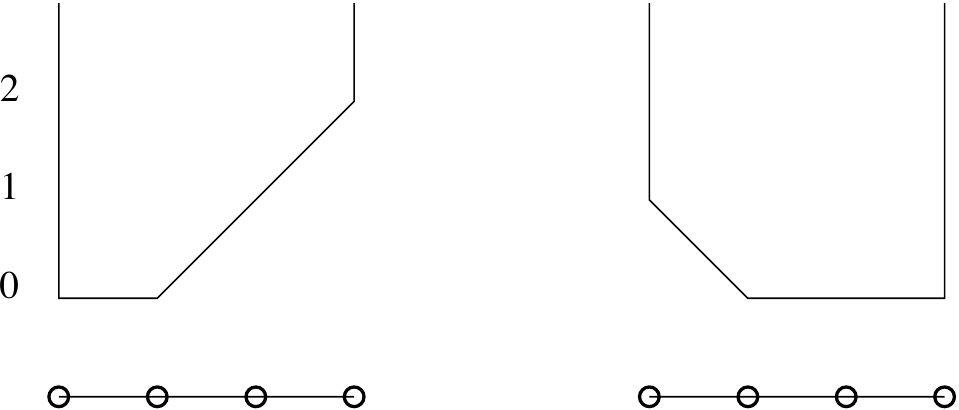}}
\caption{$\D_F$ and $\D_G$.}\label{lifts}
     \end{center}
     \end{figure}

Using the notation of ({\ref{par}}), in local coordinates the embeddings of $\tilde X_F$ and $\tilde X_G$ in $\P^3\times\C$ are $([1,t,\la t^2,\la^2 t^3],\la)$ and $([\la,t,t^2,t^3],\la)$, while in homogeneous coordinates these are 
$$
([v_1^3,v_0v_1^2,\la v_0^2v_1,\la^2 v_0^3],\la)
$$ 
and 
$$
([\la v_1^3,v_0v_1^2,v_0^2v_1,v_0^3],\la).
$$
We therefore observe that $\tilde X_F$ and $\tilde X_G$ have different parametric equations, nevertheless it is easy to see
that both of them are defined in $\P^3\times\C$ by the equations
$$
x_0x_2-\et x_1^2=0, x_1x_3-x_2^2=0, x_0x_3-\et x_1x_2=0,
$$
where $\et$ is the non-homogeneous coordinate in $\C$. These equations can also be found applying elimination theory algorithms to the two parametrizations in homogeneous coordinates, computations can be performed by hand or using computer algebra systems. 

\subsection{\bf $\P^1\times\P^1$ blown up in a point}

\begin{figure}[htbp]
     \begin{center}
     \scalebox{1}{\includegraphics{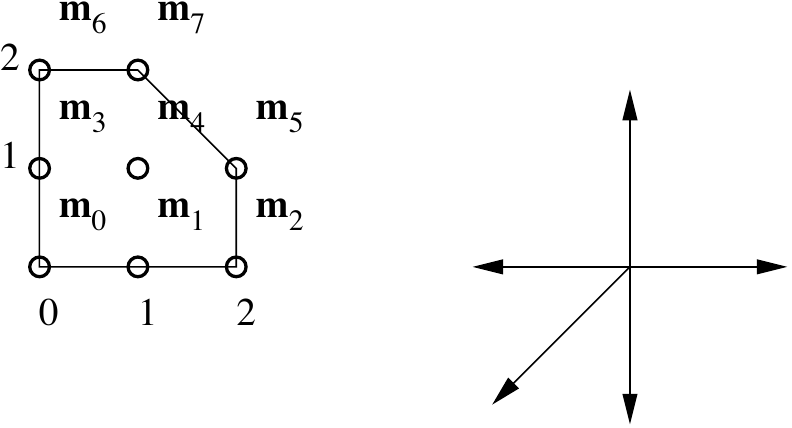}}
\caption{$\P^1\times\P^1$ blown up in a point and its normal fan.}\label{ppb}
     \end{center}
\end{figure}

Consider the polytope $\D$ in figure \ref{ppb} with its associated normal fan. 
The toric surface $X$ determined by $\D$ is $\P^1\times\P^1$ blown up in a point and embedded in $\P^7$. In local coordinates it is the closure of the image of
\begin{eqnarray}
\varphi: (\C^*)^2 & \rightarrow & \P^7 \nonumber \\
          (t_1,t_2) & \mapsto & [1,t_1,t_1^2,t_2,t_1t_2,t_1^2t_2,t_2^2,t_1t_2^2], \nonumber
\end{eqnarray}
while taking homogeneous coordinates $v_0,\ldots,v_4$ for $X_{\D}$ (one for each facet of $\D$), the embedding is

\begin{eqnarray}
X_{\D} & \hookrightarrow & \P^7 \nonumber  \\
(v_0,\ldots,v_4) & \mapsto & [ v_2^2v_3^3v_4^2, v_0v_2v_3^2v_4^2, v_0^2v_3v_4^2, v_1v_2^2v_3^2v_4 ,v_0v_1v_2v_3v_4, \label{embppb} \\
                     &         & v_0^2v_1v_4, v_1^2v_2^2v_3, v_0v_1^2v_2v_4]). \nonumber
\end{eqnarray}

Consider the semi-stable partition $\{\D_1,\D_2\}$ of $\D$:
\begin{figure}[h]
     \begin{center}
     \scalebox{1}{\includegraphics{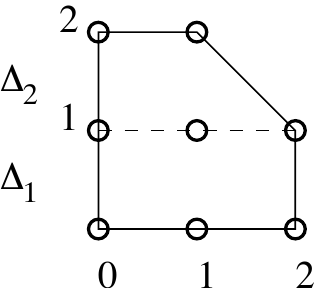}}
\caption{A semistable partition of $X$.}\label{degppb}
     \end{center}
\end{figure}

This partition gives the semi-stable degeneration of $X$ to the union of two surfaces $X_1\cup X_2$, where $X_1=\P^1\times\P^1$ and $X_2=\F^1$. 

The two possible lifting functions are
$$
F(\m_j)=
\left\{
\begin{array}{cc}
0 & j=0,\ldots,5 \\
1 & j=6,7
\end{array}
\right., \,
G(\m_j)=
\left\{
\begin{array}{cc}
1 & j=0,1,2 \\
0 & j=3,\ldots,7
\end{array}
\right..
$$
In local coordinates the embdeddings of $\tilde X_F$ and $\tilde X_G$ in $\P^7\times\C$ are 
$$
([1,t_1,t_1^2,t_2,t_1t_2,t_1^2t_2,\la t_2^2,\la t_1t_2^2],\la)
$$
 and
$$
([\la,\la t_1,\la t_1^2,t_2,t_1t_2,t_1^2t_2,t_2^2,t_1t_2^2],\la).
$$ 
We have embeddings
$$
\begin{array}{cll}
\iota_F : \tilde X_F & \hookrightarrow & \P^7\times\C  \\
(v_0,\ldots,v_4,\la) & \mapsto & ([ v_2^2v_3^3v_4^2, v_0v_2v_3^2v_4^2, v_0^2v_3v_4^2, v_1v_2^2v_3^2v_4 ,v_0v_1v_2v_3v_4,  \\
                     &         & v_0^2v_1v_4, \la v_1^2v_2^2v_3,\la v_0v_1^2v_2v_4],\la), 
\end{array}
$$
and
$$
\begin{array}{cll}
\iota_G : \tilde X_G & \hookrightarrow & \P^7\times\C   \\
(v_0,\ldots,v_4,\la) & \mapsto & ([\la v_2^2v_3^3v_4^2,\la v_0v_2v_3^2v_4^2, \la v_0^2v_3v_4^2, v_1v_2^2v_3^2v_4 ,v_0v_1v_2v_3v_4,  \\
                     &         & v_0^2v_1v_4, v_1^2v_2^2v_3, v_0v_1^2v_2v_4],\la). 
\end{array}
$$

$\tilde X_F$ and $\tilde X_G$ have different parametric equations. We find that $\tilde X_F, \tilde X_G$ are both defined in $\P^7\times\C$ by the following nine quadratic equations:
$$
\begin{array}{c}
x_3x_5-x_4^2=0,x_2x_6-\la x_4^2=0, x_1x_6-\la x_3x_4=0 \\
x_1x_5-x_2x_4=0, x_1x_4-x_2x_3=0, x_0x_6-\la x_3^2=0 \\
x_0x_5-x_2x_3=0, x_0x_4-x_1x_3=0, x_0x_2-x_1^2=0. 
\end{array}
$$
Omitting $\la$ in these equations we obtain a set of equation for $X_{\D}$ embedded in $\P^7$: these are the same equations one can compute from (\ref{embppb}) trough elimination.

\section{\bf Main results}

We use the notation of the previous sections.

\bigskip

Let $I_F$ be the ideal of all polynomials in the coordinates $x_0,\ldots,x_{\l},\et$ homogeneous in $x_0,\ldots,x_{\l}$ and vanishing on $\tilde X_F$, where $\et$ is the non-homogeneous coordinate in $\C$. In analogy with the compact case we use the notation
$$
\z^{\u}=x_0^{u_0}\ldots x_{\l}^{u_\l}\et^{u_{\l+1}},
$$
with $\u=(u_0,\ldots,u_{\l},u_{\l+1})\in\Z^{\l+2}$.

Consider the $(n+2)\times (\l+2)$  matrix
$$
\Bc^+=\Bc_F^+ =
\left(
\begin{array}{ccccc}
1        & 1      & \ldots & 1         & 0     \\
m_{10}   & m_{11} & \ldots & m_{1\l}   & 0     \\
\vdots   & \vdots &        & \vdots    &\vdots \\
m_{n0}   & m_{n1} & \ldots & m_{n\l}   & 0     \\
F(\mz)   & F(\mu) & \ldots & F(\ml)    & 1
\end{array}
\right).
$$

\begin{lemma}
$I_F$ is the linear span of all binomials $\z^{\u}-\z^{\v}$ with vectors $\u,\v\in\N^{\l+2}$ such that $\Bc^+\,\u=\Bc^+\,\v$. 
\end{lemma}

{\it Proof}. We follow Theorems 2.1 and 2.2  \cite{Sottile}.

\medskip

A binomial $\z^{\u}-\z^{\v}$, with $\u,\v\in\N^{\l+2}$, vanishing on $\psi((\C^*)^n\times\C)$ needs to be homogeneous in the coordinates $x_0,\ldots,x_{\ell}$, i.e. 
\begin{equation}
\label{hom}
\sum_{i=0}^{\ell}u_i=\sum_{i=0}^{\ell}v_i.
\end{equation}
Therefore we prove that $I_F$ is the linear span of all binomials $\z^{\u}-\z^{\v}$ with vectors $\u,\v$ such that (\ref{hom}) holds and $\Bc\u=\Bc\v$, where
$$
\Bc=\Bc_F =
\left(
\begin{array}{ccccc}
m_{10}   & m_{11} & \ldots & m_{1\ell} & 0 \\
\vdots   & \vdots &        & \vdots    & \vdots \\ 
m_{n0}   & m_{n1} & \ldots & m_{n\ell} & 0 \\
F(\mz)   & F(\mu) & \ldots & F(\ml)    & 1
\end{array}
\right).
$$

\smallskip

Consider a monomial $\z^{\u}$ and restrict it to $\psi((\C^*)^n\times\C)$:
\begin{eqnarray}
 \z^{\u}_{|\psi((\C^*)^n\times\C)} &=& (x_0^{u_0}\ldots x_{\l}^{u_{\l}}\et^{u_{\l+1}})_{|\psi((\C^*)^n\times\C)}=   \nonumber \\
                                  & = & (t_1^{m_{10}}\ldots t_n^{m_{n0}}\la^{F(\m_0)})^{u_0}\ldots(t_1^{m_{1\l}}\ldots t_n^{m_{n\l}} \la^{F(\m_{\l})})^{u_{\l}}\cdot \nonumber \\ 
& \cdot &\la^{u_{\l+1}}=  \nonumber \\  
                                  & = & t_1^{m_{10}u_0+\ldots+m_{1\l}u_{\l}}\ldots t_n^{m_{n0}u_0+\ldots+ m_{n\l}u_{\l}}\cdot \nonumber\\
 &  \cdot &  \la^{F(\m_0)u_0+\ldots+ F(\m_{\l})u_{\l}+u_{\l+1}}=  \nonumber \\
                                  & = &  T^{\Bc\u},   \nonumber
\end{eqnarray}
with $T=(t_1,\ldots,t_n,\la)$. 

This shows that in the hypothesis (\ref{hom}), $\z^{\u}-\z^{\v}$ vanishes on $\psi((\C^*)^n\times\C)$ (and hence belongs to $I_F$) if and only if $\Bc\u=\Bc\v$.

\bigskip

Now we show that these binomials generate $I_F$ as a $\C$-vector space: we follow Strumfels's book \cite{StrumL}. Strumfels considers the (compact) toric variety defined as in (\ref{toric}) and doesn't deal with the homogeneous vs. non-homogeneous question.

\bigskip

Fix a monomial ordering $>$ on $\C[x_0,\ldots,x_{\l},\et]$, and remember that this is a well-ordering on the set of monomials $\z^{\u}$. Suppose the set $R$ of polynomials $f\in I_F$ which cannot be written as a $\C$-linear combination of binomials as above is non-empty and take $f\in R$ such that
$$
\mbox{\emph{LM}}_{>}(f)=\min_{g\in R}\mbox{\emph{LM}}_{>}(g),
$$ 
where $\mbox{\emph{LM}}_{>}(f)$ is the leading monomial of $f$ with respect to $>$.
We can suppose $f$ to be monic, so that its leading term $\mbox{\emph{LT}}_{>}(f)$ is its leading monomial, let this be the monomial $\z^{\u}$. 

When we restrict $f$ to $\psi((\C^*)^n\times\C)$ we get an expression containing $T^{\Bc\u}$ as a term and which is equal to zero. Hence the term $T^{\Bc\u}$ must cancel in this expression. This means that there is some other monomial $\z^{\v}$ appearing in $f$ such that $\Bc \u=\Bc \v$ and (\ref{hom}) holds.

Moreover $\z^{\u}>\z^{\v}$. The polynomial 
$$
f':=f-\z^{\u}+\z^{\v}
$$
belongs to $I_F$ and to $R$ but since $\mbox{\emph{LM}}_>(f)>\mbox{\emph{LM}}_>(f')$, we get a contradiction. $\quad\Box$

\begin{theorem}
\label{ker}
$I_F=\langle \z^{\u^+}-\z^{\u^-}| \u\in\ker(\Bc^+) \mbox{ \emph{and} } \u\in\Z^{\l+2}\rangle.$
\end{theorem}

{\it Proof}. On one hand, $\u\in\ker(\Bc^+)$ if and only if $\Bc^+ \u^+=\Bc^+ \u^-$. On the other hand we show that if $\Bc^+ \v=\Bc^+ \w$ (and (\ref{hom}) holds), then $\z^\v-\z^\w=h(\z^{\u^+}-\z^{\u^-})$, for some polynomial $h$ and vector $\u\in\ker(\Bc^+)\cap\Z^{\l+2}$; the statement will then follow from the theorem. 

If $\Bc^+ \v=\Bc^+ \w$, then  $\v-\w\in\ker(\Bc^+)$.
\begin{eqnarray}
\z^\v-\z^\w=\z^\w(\z^{\v-\w}-1) & = & \z^\w\z^{-(\v-\w)^-}(\z^{(\v-\w)^+} - \z^{(\v-\w)^-}) \nonumber \\
                           & = & \z^{\w-(\v-\w)^-}(\z^{(\v-\w)^+}-\z^{(\v-\w)^-}) \nonumber
\end{eqnarray}
It is easy to show that $\w-(\v-\w)^-\in\N^{\l+2}$. $\quad\Box$ 

\bigskip

Now let $G$ be the second lift, then we can consider the matrix
$\Bc^+_G$ and characterize the toric ideal $I_G$ of $\tilde X_G$ as
above. In general $\tilde X_G$ will have a different parametrization
from the one of $\tilde X_F$, moreover the normal fans are
different.

\bigskip

Our main result is

\bigskip

\begin{theorem}
$\tilde X_F$ and $\tilde X_G$ have the same equations in $\P^\l\times\C$, i.e. $I_F=I_G$.
\end{theorem}

{\it Proof}. Reorder the $\m_j$'s such that $\mz,\ldots,\m_r\in\D_1-\D_2$, $\m_{r+1},\ldots,\m_s\in\D_1\cap\D_2$ and $\m_{s+1},\ldots,\m_\l\in\D_2-\D_1$, then we have
$$
\Bc_F^+=
\left(
\begin{array}{cccccccccc}
1      &..& 1      & 1        &..& 1      & 1                 &..& 1             & 0 \\
m_{10} &..& m_{1r} & m_{1,r+1}&..& m_{1s} & m_{1,s+1}         &..& m_{1\l}       & 0 \\  
\vdots &  & \vdots & \vdots   &  & \vdots & \vdots            &  & \vdots        & \vdots \\ 
m_{n0} &..& m_{nr} & m_{n,r+1}&..& m_{ns} & m_{n,s+1}         &..& m_{n\l}       & 0 \\
0      &..& 0      & 0        &..& 0      & L_F(\m_{s+1})     &..& L_F(\m{\l})   & 1
\end{array}
\right),
$$
and  
$$
\Bc_G^+=
\left(
\begin{array}{cccccccccc}
1        &..& 1         & 1         &..& 1      & 1         &..& 1       & 0 \\
m_{10}   &..& m_{1r}    & m_{1,r+1} &..& m_{1s} & m_{1,s+1} &..& m_{1\l} & 0 \\ 
\vdots   &  & \vdots    & \vdots    &  & \vdots & \vdots    &  & \vdots  & \vdots \\
m_{n0}   &..& m_{nr}    & m_{n,r+1} &..& m_{ns} & m_{n,s+1} &.. & m_{n\l} & 0 \\
L_G(\m_0)&..& L_G(\m_r) & 0         &..& 0      & 0         &..& 0       & 1
\end{array}
\right).
$$

Let $E$ be the $(n+2)\times (n+2)$ elementary matrix 
$$
\left(
\begin{array}{ccccc}
1       & 0       & \ldots & 0      & 0 \\   
0       & 1       & \ldots & 0      & 0 \\
\vdots  & \vdots  &        & \vdots & \vdots  \\
0       & 0       & \ldots & 1      & 0 \\
a_{n+1} & a_1     & \ldots & a_n    & 1 \\
\end{array}
\right)
\in SL_{n+2}(\Z)
$$ 

we have 
$$
E\cdot\Bc_G^+=\Bc_F^+,
$$
and hence 
$$
\ker\Bc_F^+=\ker\Bc_G^+.
$$
The theorem follows from Theorem (\ref{ker}). $\quad\Box$

\bigskip

Going back to the examples above, if $X$ is the twisted cubic, we have
$$
\Bc_F^+=
\left(
\begin{array}{ccccc}
1 & 1 & 1 & 1 & 0 \\
0 & 1 & 2 & 3 & 0 \\
1 & 0 & 0 & 0 & 1
\end{array}
\right),\,
\Bc_G^+=
\left(
\begin{array}{ccccc}
1 & 1 & 1 & 1 & 0 \\
0 & 1 & 2 & 3 & 0 \\
0 & 0 & 1 & 2 & 1
\end{array}
\right),
$$
and $E$ is the $3\times 3$ elementary matrix
$$
\left(
\begin{array}{ccc}
1 & 0 & 0  \\
0 & 1 & 0  \\
1 & -1 & 1
\end{array}
\right)\in SL_3(\Z).
$$

\bigskip

In the case of $\P^1\times\P^1$ blown up in a point, we have
$$
\Bc_F^+=
\left(
\begin{array}{ccccccccc}
1 & 1 & 1 & 1 & 1 & 1 & 1 & 1 & 0 \\
0 & 1 & 2 & 0 & 1 & 2 & 0 & 1 & 0 \\
0 & 0 & 0 & 1 & 1 & 1 & 2 & 2 & 0 \\
1 & 1 & 1 & 0 & 0 & 0 & 0 & 0 & 1
\end{array}
\right),
$$ 
$$
\Bc_G^+=
\left(
\begin{array}{ccccccccc}
1 & 1 & 1 & 1 & 1 & 1 & 1 & 1 & 0 \\
0 & 1 & 2 & 0 & 1 & 2 & 0 & 1 & 0 \\
0 & 0 & 0 & 1 & 1 & 1 & 2 & 2 & 0 \\
0 & 0 & 0 & 0 & 0 & 0 & 1 & 1 & 1
\end{array}
\right),
$$
and 
$$
E=
\left(
\begin{array}{cccc}
1 & 0 & 0 & 0 \\
0 & 1 & 0 & 0 \\
0 & 0 & 1 & 0 \\
1 & 0 & -1 & 1
\end{array}
\right)\in SL_4(\Z).
$$

It would be interesting to extend such results to semi-stable partitions of a polytope $\D$ in an arbitrary number of subpolytopes.

\newpage

\bigskip
\noindent Marina Marchisio\\
Universit\`a di Torino\\
Dipartimento di Matematica\\
Via Carlo Alberto, 10\\
10123 Torino (Italy) \\
marina.marchisio@unito.it

\bigskip
\noindent Vittorio Perduca\\
Universit\`a di Torino\\
Dipartimento di Matematica\\
Via Carlo Alberto, 10\\
10123 Torino (Italy) \\
vittorio.perduca@unito.it

\end{document}